%%%%%%%%%%%%%%%%%%%%%
%%This contains amstex file for the survey paper 
%%"Surgery and stratified spaces"
%%by Hughes and Weinberger for the Wall book.
%% This was begun on 9 Dec 1996
%% Submitted version 2 July 1997
%% On July 29, 1997 I enlarged the mag. and corrected Fillmore/Rothenberg
%% in bibliography.
\input amstex
\documentstyle{amsppt}
\NoBlackBoxes
\magnification 1100
%\input symb.tex
%%%****************************************************************
%%%%%%%%%%%%%%%%%%%%%%%%%%%%%%%%%%%%%%%%%%%%%%%%%%%%%%%%%%%%%%%%%%%%%%%%
%%%%%%%%%%%%%%%%%%%%%

\define\bq{\Bbb Q}
\define\br{\Bbb R}
\define\bz{\Bbb Z}
\define\cl{\operatorname{cl}}
\define\ov{\overline}
\define\cyl{\operatorname{cyl}}
\define\ho{\operatorname{holink}}
\define\inr{\operatorname{int}}
\define\rel{\text{\rm rel\ }}
\define\oc{\overset\circ\to{\text{\rm c}}}
\define\id{\operatorname{id}}
\define\Wh{\operatorname{Wh}}

\topmatter
\title  Surgery and Stratified Spaces\endtitle
%\rightheadtext{}
\author Bruce Hughes and Shmuel Weinberger\endauthor
\address Department of Mathematics, Vanderbilt University, Nashville,
Tennessee 37240\endaddress
\email hughes\@math.vanderbilt.edu\endemail
\address Department of Mathematics, University of Chicago, 
Chicago, Illinois 60637\endaddress
\email shmuel\@math.uchicago.edu\endemail

%\abstract  \endabstract
\thanks The first author was 
supported in part by NSF Grant DMS--9504759\endthanks
\thanks The second author was supported in part by NSF Grant DMS--9504913
\endthanks
%\date 2 july 1997\enddate
\endtopmatter

\document
\head 0. Introduction
\endhead
The past couple of decades has seen significant progress in the
theory of stratified spaces through the application of controlled
methods as well as through the applications of intersection
homology.  In this paper we will give a cursory introduction to this
material, hopefully whetting your appetite to peruse more thorough
accounts.

        In more detail, the contents of this paper are as follows:  the
first section deals with some examples of stratified spaces and
describes some of the different categories that have been considered
by various authors.  For the purposes of this paper, we will work in
either the PL category or a very natural topological category
introduced by Quinn \cite{Q4}.  The next section discusses intersection
homology and how it provides one with a rich collection of self dual
sheaves.  These can be manipulated by ideas long familiar to surgery
theorists who have exploited Poincar\'e duality from the start.  We
will give a few applications of the tight connection between an
important class of stratified spaces (Witt spaces), self dual
sheaves, and K-theory; one last application will appear in the final
section of the paper (where we deal with the classification of
``supernormal'' spaces with only even codimensional strata).

        Section three begins an independent direction, more purely
geometric.  We describe the local structure of topological stratified
spaces in some detail, in particular explaining the teardrop
neighborhood theorem (\cite{HTWW}, \cite{H2}) and giving applications to
isotopy theorems and the like.  The last three sections describe the
theory of surgery on stratified spaces, building on our understanding
of teardrop neighborhoods, and some applications to classification
problems (other applications can also be found in the survey 
\cite{CW4}).

\head 1. Definitions and Examples of Stratified Spaces 
\endhead
A {\it stratification\/} $\Sigma =\{X_i\}$ of a space $X$ is a locally
finite decomposition of $X$ into pairwise disjoint, locally closed
subsets of $X$ such that each $X_i\in\Sigma$ is a topological manifold.
We always assume that $X$ is a locally compact, separable metric space 
and that $\Sigma$ satisfies the {\it Frontier Condition\/}:
$\cl X_i\cap X_j\not=\emptyset$ if and only if
$X_j\subseteq\cl X_i$. The index set is then partially ordered by
$j\leq i$ if and only if $X_j\subseteq\cl X_i$. The set $X_i\in\Sigma$
is called a {\it stratum\/} and $X^i =\cl X_i =\cup\{ X_j ~|~ j\leq i\}$ is a
{\it skeleton\/} (or {\it closed stratum\/} in the terminology of \cite{W4}).

Partitioning non-manifold spaces into manifold pieces is a very old idea ---
one has only to consider polyhedra in which the strata are differences
between successive skeleta. However, it was not until relatively recently
that  attention was paid to how
the strata should fit together, or to
the geometry of the neighborhoods of strata. 
In 1962  Thom \cite{T1} discussed
stratifications in which the strata have neighborhoods
which fibre over that stratum and which have ``tapis'' maps (the
precursor to the tubular maps in Mather's formulation in 1.2 below).
It was also in this paper that Thom conjectured that the 
topologically stable maps between two smooth manifolds are dense in the
space of all smooth maps. In fact, it was Thom's program for attacking that
conjecture which led him to a study of stratifications \cite{T2}. The connection
between stratifications and topological stability
(and, more generally, the theory of singularities of smooth maps) is
outside the scope of this paper, but the connections have continued 
to develop
(for a recent account, see the book of du Plessis and Wall \cite{dPW}.)

Here we review the major formulations of the conditions on neighborhoods
of strata. These are due to Whitney, Mather, Browder and Quinn,
Siebenmann, and Quinn. The approaches of Whitney, Mather, Browder and
Quinn are closely related to Thom's original ideas. These types of 
stratifications are referred to as {\it geometric stratifications\/}.
The approaches of Siebenmann and Quinn are attempts at finding an
appropriate topological
setting.

\subhead 1.1 Whitney stratifications\endsubhead
In two fundamental papers \cite{Wh1},\cite{Wh2},
Whitney clarified some of Thom's ideas on
stratifications and introduced his Conditions A and B.
To motivate these conditions consider a real algebraic set
$V\subseteq\br^n$, the common locus of finitely many real polynomials.
The singular set $\Sigma V$ of all points where $V$ fails to be
a smooth manifold is also an algebraic set.
There is then a finite filtration 
$V=V^m\supseteq V^{m-1}\supseteq \cdots\supseteq V^0\supseteq V^{-1}=\emptyset$
with $V^{i-1}=\Sigma V^i$. One obtains a stratification of $V$
by considering the strata $V_i=V^i\setminus V^{i-1}$.
However, with this naive construction the strata need not have 
geometrically well-behaved neighborhoods; that is, the local topological
type need not be locally constant along strata. For example, consider
the famous Whitney umbrella
which is the locus of $x^2=zy^2$, an algebraic set in $\br^3$.
The singular set $\Sigma V$ is the $z$-axis and is a smooth manifold,
so one obtains just two strata, $V\setminus\Sigma V$ and $\Sigma V$.
However, there is a drastic change in the neighborhood of $\Sigma V$ in
$V$ as one passes through the origin: for negative $z$ small neighborhoods
meet only $\Sigma V$ whereas this is not the case for positive $z$.

If $X$, $Y$ are smooth submanifolds of a smooth manifold $M$, then
$X$ is {\it Whitney regular over\/} $Y$ if whenever $x_i\in X$,
$y_i\in Y$ are sequence of points converging to some $y\in Y$, the lines
$l_i=\ov{x_iy_i}$ converge to a line $l$, and the tangent spaces
$T_{x_i}X$ converge to a space $\tau$, then
\roster
\item"(A)" $T_yY\subseteq\tau$ and
\item"(B)" $l\subseteq\tau$.
\endroster
A stratification $\Sigma=\{X_i\}$ of $X$ is a {\it Whitney stratification\/}
if whenever $j\leq i$, $X_i$ is Whitney regular over $X_j$.

In the Whitney umbrella $V$, $V\setminus\Sigma V$ is not
Whitney regular over $\Sigma V$ at the origin.
However, the stratification can be modified to give a Whitney stratification
and a similar construction works for a class of spaces 
more general than algebraic sets:
a subset $V\subseteq\br^n$ is a {\it semi-algebraic set\/} if it is a 
finite union of sets which are the common solutions of finitely many
polynomial equations and inequalities. Examples include real algebraic sets
and polyhedra. In fact, the class of semi-algebraic sets is the smallest class 
of euclidean subsets 
containing the real algebraic sets and which is  closed under 
images of linear projections. 
If $V$ is semi-algebraic, then there is a finite filtration as in the case
of an algebraic set discussed above obtained by considering iterated
singular sets. This filtration can be modified by removing from the strata
the closure of the set of points where the Whitney conditions fail to hold.
In this way, semi-algebraic sets are given Whitney stratifications
(see \cite{GWdPL}). 

In fact, Whitney \cite{Wh2} showed that any real or complex analytic set
admits a Whitney stratification.  These are sets defined analogously to
algebraic sets with analytic functions  used instead of polynomials.
Lojasiewicz \cite{Lo} then showed that semi-analytic sets (the analytic 
analogue
of semi-algebraic sets) have Whitney
stratifications. An even more general class of spaces, namely the 
{\it subanalytic sets\/}, were shown by Hardt \cite{Hr} to admit Whitney
stratifications. For a modern and thorough discussion of stratifications for
semi-algebraic and subanalytic sets see Shiota \cite{Shi}.

\subhead 1.2 Mather stratifications: tube systems\endsubhead
Mather clarified many of the ideas of Thom and Whitney and gave complete
proofs of the isotopy lemmas of Thom. He worked with a definition of
stratifications closer to Thom's original ideas than to Whitney's, but then
proved that spaces with Whitney stratifications are stratified in his sense.

\definition{Definition} For 
$0\leq k\leq +\infty$,
a {\it Mather $C^k$-stratification\/} of $X$ is  a triple $(X,\Sigma,${\bf T}$)$
such the following hold:
\roster
\item $\Sigma =\{X_i\}$ is a stratification of $X$ 
such that  each stratum $X_i\in\Sigma$ is a $C^k$-manifold.
\item {\bf T} $=\{T_i,\pi_i,\rho_i\}$ is called a {\it tube system\/}
and $T_i$ is an open neighborhood of $X_i$ in $X$, called  the
{\it tubular neighborhood\/} of $X_i$,
$\pi_i:T_i\to X_i$ is  a retraction, called the {\it local retraction\/}
of $T_i$, and $\rho_i:T_i\to [0,\infty)$ is a map
such that $\rho_i^{-1}(0) = T_i$.

\item For each $X_i, X_j\in\Sigma$, if $T_{ij} =T_i\cap X_j$ and
the restrictions of $\pi_i$ and $\rho_i$ to $T_{ij}$ are denoted
$\pi_{ij}$ and $\rho_{ij}$, respectively, then
the map
$$(\pi_{ij},\rho_{ij}):T_{ij}\to X_i\times (0,\infty)$$
is a $C^k$-submersion.

\item If $X_i,X_j,X_k\in\Sigma$, then the following compatibility 
conditions hold for
$x\in T_{jk}\cap T_{ik}\cap \pi_{jk}^{-1}(T_{ij})$:
$$\pi_{ij}\circ\pi_{jk}(x)=\pi_{ik}(x),$$
$$\rho_{ij}\circ\pi_{jk}(x)=\rho_{ik}(x).$$
\endroster
\enddefinition

When $k=0$ above, a $C^0$-submersion, or topological
submersion, means every point in the domain has a neighborhood in
which the map is topologically equivalent to a projection (see \cite{S2}).

Mather \cite{Ma1}, \cite{Ma2} proved that Whitney stratified spaces
have Mather $C^\infty$-stratifications.

The Thom isotopy lemmas mentioned above are closely related to the geometric
structure of neighborhoods of strata.
For example, the first isotopy lemma says that if $f:X\to Y$ is a proper
map between Whitney stratified spaces with the property that for every
stratum  $X_i$ of $X$ there exists a stratum $Y_j$ of $Y$ such that
$f|:X_i\to Y_j$ is a smooth submersion, then $f$ is a fibre bundle
projection (topologically --- not smoothly!) and has local trivializations
which preserve the strata.  Mather applied this to the {\it tubular maps\/}
$$\pi_i\times\rho_i:T_i\to X_i\times [0,\infty)$$
defined on the tubular neighborhoods of the strata of a Whitney stratified
space $X$ in order to show that every stratum $X_i$  has a neighborhood
$N$ such that the pair $(N,X_i)$ is homeomorphic to $(\cyl(f),X_i)$
where $\cyl(f)$ is the mapping cylinder of some fibre bundle projection
$f:A\to X_i$. In fact, $\cl X_i$ has a mapping cylinder neighborhood
given by a map which is a smooth submersion over each stratum.
The existence of these mapping cylinder neighborhoods
was abstracted by Browder and Quinn as is seen next.

\subhead 1.3 Browder-Quinn stratifications: mapping cylinder neighborhoods
\endsubhead
In order to classify stratified spaces Browder and Quinn \cite{BQ}
isolated the mapping cylinder structure as formulated by Mather.
The mapping cylinder  was then part of the data that was to
be classified. More will be said about this in \S4 below. Here we
recall their definition.

Let $\Sigma=\{ X_i\}$ be a stratification of a space $X$ such that each
stratum $X_i$ is a $C^k$-manifold. The {\it singular set\/} $\Sigma X_i$
is $\cl X_i\setminus X_i=\cup\{ X_j ~|~ j < i\}$. (This is in general
bigger than the singular set as defined in 1.1.)

\definition{Definition} 
$\Sigma$ is a $C^k$ {\it geometric stratification\/} of $X$ if for every 
$i$ there is a closed neighborhood $N_i$ of $\Sigma X_i$ in $X^i=\cl X_i$ and a
map $\nu_i:\partial N_i\to\Sigma X_i$ such that
\roster
\item $\partial N_i$ is a codimension $0$ submanifold of $X_i$,
\item $N_i$ is the mapping cylinder of $\nu_i$ (with $\partial N_i$ and
$\Sigma N_i$ corresponding to the top and bottom of the cylinder),
\item if $j< i$ and $W_j = X_j\setminus \inr N_j$, then
$\nu_i|:\nu_i^{-1}(W_j)\to W_j$ is a $C^k$-submersion.
\endroster
\enddefinition

The complement of $\inr N_i$ in $X^i$ is called a {\it closed pure
stratum\/} and is denote $\ov X^i$.

Note this definition incorporates a topological theory by taking $k=0$.
Browder and Quinn also pointed out that by relaxing the condition
on the maps $v_i$ other theories can be considered. For example, one
can insist that the
strata be PL manifolds and the $v_i$ be PL block bundles with manifold
fibers.

\subhead 1.4 Siebenmann stratifications: local cones\endsubhead
In the late 1960s Cernavski \cite{Ce} developed intricate geometric techniques
for deforming homeomorphisms of topological manifolds. In particular,
he proved that the group of self homeomorphisms of a compact manifold
is locally contractible by showing that two sufficiently
close homeomorphisms are canonically isotopic.
The result was reproved by Edwards and Kirby \cite{EK}
by use of Kirby's torus trick. Siebenmann \cite{S2} developed the technique
further in order to establish the local contractibility of homeomorphism
groups for certain nonmanifolds, especially, compact polyhedra.

Siebenmann's technique applied most naturally to stratified spaces and
a secondary aim of \cite{S2} was to introduce a class of stratified
spaces that he thought might
``come to be the topological analogues of polyhedra in the 
piecewise-linear realm or of Thom's stratified sets in the
differentiable realm.'' These are the {\it locally conelike TOP
stratified sets\/} whose defining property is that strata 
are topological manifolds and for each point
$x$ in the $n$--stratum there is a compact locally conelike TOP
stratified set $L$ (with fewer strata --- the definition is inductive)
and a stratum-preserving homeomorphism of $\br^n\times\oc L$ onto
an open neighborhood of $x$ where $\oc L$ denotes the open cone on $L$
and the homeomorphism takes $0\times\{\text{vertex}\}$ to $x$.
Simple examples include polyhedra and the topological ($C^0$) versions

It is important to realize that Siebenmann didn't just take the topological
version of Mather's stratified space, but he did have Mather's
$C^0$-tubular maps {\it locally\/} at each point. The reason he was able
to work in this generality was that the techniques for proving
local contractibility of homeomorphism groups were purely local. 

As an example, consider a pair $(M,N)$ of topological manifolds
with $N$ a locally flat submanifold of $M$. With the two strata,
$N$ and $M\setminus N$, the local flatness verifies that this is a 
locally conelike stratification. However, Rourke and Sanderson \cite{RS}
showed that $N$ need not have a neighborhood given by the mapping cylinder
of a fibre bundle projection. Thus, Siebenmann's class is definitely
larger than the topological version of the Thom-Whitney-Mather
class.

On the other hand, Edwards \cite{E} did establish that locally flat 
submanifolds of high dimensional topological manifolds do, in general,
have mapping cylinder neighborhoods. However, the maps to the submanifold
giving the mapping cylinder need not be a fibre bundle projection.
It turns out that the map is a {\it manifold approximate fibration\/}, 
a type of map which figures prominently in the discussion of the geometry
of homotopically stratified spaces below.

Later, Quinn \cite{Q2,II} and Steinberger and West \cite{StW} gave examples of
locally conelike TOP stratified sets in which the strata do no have 
mapping cylinder neighborhoods of any kind. In fact, their examples are
orbit spaces of finite groups acting locally linearly on topological
manifolds. Such orbit spaces are an important source of examples of
locally conelike stratified sets and many of advances in the theory of 
stratified spaces were made with applications to locally linear actions 
in mind. These examples were preceded by an example mentioned by Siebenmann
\cite{S3} of a locally triangulable non-triangulable space.

Milnor's counterexamples to the Hauptvermutung \cite{M1} give non-homeomorphic
polyhedra whose open cones are homeomorphic. As Siebenmann observed, these
show that the links in
locally conelike stratified sets are not unique up to homeomorphism. 
Siebenmann does prove that the links are stably homeomorphic after crossing
with a euclidean space of the dimension of the stratum plus $1$. 
The non-uniqueness
of links points to the fact that Siebenmann's stratified spaces are too rigid
to really be the topological analogue of polyhedra and smoothly stratified 
sets, whereas the stable uniqueness foreshadows the uniqueness up to
controlled homeomorphism of fibre germs of manifold approximate 
fibrations \cite{HTW1}.

The main applications obtained by Siebenmann, namely local contractibility
of homeomorphism groups, isotopy extension theorems, and the fact
that many proper submersions are bundle projections, can all be 
generalized to the setting of homotopically stratified sets discussed below.

Siebenmann himself experimented with a less rigid class of stratified spaces,
called {\it locally weakly conelike\/}. In order to include in this class
stratified spaces with isolated singularities which arise as 
the one-point compactifications of manifolds with nonvanishing
Siebenmann obstruction \cite{S},
he no longer required the existence of links.
However, neighborhoods around strata
of dimension $n$ were still required to split off a factor of $\br^n$
locally.
In other words, in a locally conelike set $X$ a point in the $n$-dimensional
stratum $X_n$ has a neighborhood $U$ in $X$ with $U\setminus X_n$ homeomorphic
to $L\times\br^{n+1}$ with $L$ the compact link. 
In a weakly conelike set, $U\setminus X_n$ is
homeomorphic to $C\times\br^n$ with $C$ a noncompact space with a certain
tameness property at infinity. While this generalization was a move in the 
right direction, the role of the weak link $C$ prevented further 
developments and it was left to Quinn to make a bolder generalization.

\subhead 1.5 Quinn stratifications: homotopy mapping cylinders\endsubhead
Quinn \cite{Q5} introduced a class of spaces which we will call
{\it manifold homotopically stratified sets\/}. His objective was to
``give a setting for the study of purely topological stratified
phenomena, particularly group actions on manifolds.'' As has been
pointed out above, the previously defined topologically stratified
spaces were inadequate. On one hand, the geometrically stratified
spaces (that is, the topological version
of Thom's spaces as formulated by Mather or Browder and Quinn)
require too strong of a condition on neighborhoods of strata
(namely, mapping cylinder neighborhoods) ruling out important examples
(like locally flat submanifolds and orbit spaces of locally linear
group actions). On the other hand, the locally conelike stratified
sets of Siebenmann require a very strong local condition which need
not propagate to the entire neighborhood of the strata. Without an
understanding of the geometry of neighborhoods of strata, topological
stratified versions of surgery, transversality, and $h$-cobordism theorems
were missing.

Quinn formulated his definition to be equivalent to saying that for
$j< i$, $X_i\cup X_j$ is homotopy equivalent near $X_j$ to the mapping 
cylinder of some fibration over $X_j$. This has two pleasant properties.
First, besides the geometric condition that the strata be manifolds, the
definition is giving in homotopy theoretic terms. Second, the condition
concerns neighborhoods of strata rather than closed strata, so that, in
particular, there are no complicated compatibility conditions where
more than two strata meet. The links are now defined only up to homotopy.

Even without a geometric condition on neighborhoods of strata,
Quinn was able to derive geometric results. These will be discussed in 
\S3 below along with a theorem of Hughes, Taylor, Weinberger and Williams
which says that neighborhoods of strata do carry a weak geometric structure.
One thing that Quinn did not do was to develop a surgery theory for these
manifold homotopically stratified sets. That piece of the picture
was completed by
Weinberger \cite{W4} (see \S5 below).

\subhead 1.6 Group actions\endsubhead
Suppose that $G$ is a finite group acting on a topological manifold $M$.
One attempts to study the action by studying the orbit space $M/G$ and
the map $M\to M/G$. For example, if $G$ acts freely, then $M/G$ is a manifold
and $M\to M/G$ is a covering projection. Moreover, the surgery theoretic 
set of equivariant
manifold structures on $M$ is in $1-1$ correspondence with the set of 
manifold structures on $M/G$ via the pull-back construction. 

When the action is not free, $M/G$ must be viewed as  a space with 
singularities and $M\to M/G$ as a collection of covering projections.
The prototypical example occurs when $M$ is a closed Riemann surface and
$G$ is a finite cyclic group acting analytically. Then $M\to M/G$ is a
{\it branched covering\/}.

More generally, if $M$ is a smooth manifold and $G$ acts differentiably,
then $M$ has a Whitney stratification with the strata $M_{(H)}$ indexed
by conjugacy classes of subgroups of $G$ and consisting of all points
with isotropy subgroup conjugate to $H$. This induces a Whitney
stratification of $M/G$. The standard reference is Lellmann \cite{Le}, but
Dovermann and Schultz \cite{DS} provide a more accessible proof.
In the more general setting of a compact Lie group $G$, 
Davis \cite{Dv1} showed how to view $M\to M/G$ as a collection of fibre
bundle projections based on the fact that each $M_{(H)}\to M_{(H)}/G$
is a smooth fibre bundle projection with fibre $G/H$.

Now if the action of the finite group $G$ on the topological manifold
$M$ is locally linear (also called locally smooth), then the examples
of Quinn and Steinberger and West show (as mentioned
above) that $M/G$ need not have a geometric stratification, but it is a locally
conelike TOP stratified set, and so Siebenmann's results can be applied.
Lashof and Rothenberg \cite{LR} used stratification theory of the orbit
space to classify equivariant smoothings of locally smooth $G$-manifolds.
Hsiang and Pardon \cite{HsP} and Madsen and Rothenberg \cite{MR} used
stratifications for the classification of linear representations up to
homeomorphism (see also \cite{CSSW}, \cite{CSSWW}, \cite{HP}).
Stratifications also played an important role in the work of
Steinberger and West \cite{StW} on equivariant $s$-cobordism theorems and
equivariant finiteness obstructions.

The stratification theory of the orbit space actually corresponds with
the isovariant, rather than the equivariant, theory of the manifold.

Locally linear actions on topological manifolds have the property that
fixed sets are locally flat submanifolds. It is natural to consider
all such actions. These are essentially the actions whose orbit space is a
manifold homotopically stratified set. After being introduced by
Quinn \cite{Q5}, Yan \cite{Y} applied Weinberger's stratified
surgery (see \S5 below) to study equivariant periodicity.
More recently, Beshears \cite{Bs} made precise the properties of the
map $M\to M/G$ and proved that the isovariant structures on $M$
are in $1-1$ correspondence with the stratum preserving structures
on $M/G$.

\subhead 1.7 Mapping cylinders\endsubhead
Mapping cylinders provide examples of spaces with singularities.
The mapping cylinder $\cyl(p)$ of a map $p:M\to N$ between manifolds
has a natural stratification with three strata: the top $M$, the bottom
$N$ and the space in between $M\times (0,1)$. The properties of the 
stratification depend on the map $p$.
With this stratification $\cyl(p)$ is geometrically stratified if and only
if $p\times\id_\br$ can be approximated arbitrarily closely by
fibre fibre bundle projections. On the other hand, $\cyl(p)$ is 
a manifold homotopically stratified set if and only if $p$ is a manifold
approximate fibration. The cylinder is nonsingular; i.e., a manifold
with $N$ a locally flat submanifold if and only if $p$ is a manifold
approximate fibration with spherical homotopy fibre. (Here and elsewhere
in this section, we ignore problems with low dimensional strata.)

More generally, one can consider the mapping cylinder of a map
$p:X\to Y$ between stratified spaces which take each stratum of $X$ into
a stratum of $Y$. The natural collection of strata of $\cyl(p)$ contains
the strata of $X$ and $Y$. Cappell and Shaneson \cite{ChS4} observed that
even if one considers maps between smoothly stratified spaces which are
smooth submersions over each stratum of $X$, then
$\cyl(p)$ need not be smoothly stratified (they refer to an example of Thom
\cite{T1}). However, Cappell and Shaneson \cite{CS5}
proved that such cylinders are 
manifold homotopically stratified sets, showing that the stratifications
of Quinn arise naturally in the theory of smoothly stratified spaces.

Even more generally, the mapping cylinder $\cyl(p)$ of a stratum preserving
map between manifold homotopically stratified sets is itself a
manifold homotopically stratified set 
(with the natural stratification) if and only if
$p$ is a manifold stratified approximate fibration \cite{H2}.

\comment
\subhead 1.8 Stratified Morse theory\endsubhead
The extension by Goresky and MacPherson  \cite{GM4} of classical Morse
theory for manifolds to singular spaces requires that the singular
spaces have Whitney stratifications. The extension of the
notion of a Morse function to the singular case requires consideration of
a sort of nondegeneracy condition normal to strata; hence, there must
be geometric structure on neighborhoods of strata.

Bestvina and Brady \cite{BB} have developed a Morse theory for singular
spaces which are cell complexes with compatible affine structure on each
cell. This is similar to what a PL version of Goresky and MacPherson's
theory should look like. Their motivation was to analyze homological
properties of groups by associating to the groups certain Morse
functions on affine cell complexes. Recently, Davis \cite{Dv3}
has used these techniques to construct a finitely generated Poincar\'e
duality group which is not finitely presented.
\endcomment

\head 2.  Intersection Homology and Surgery Theory\endhead
In the mid 70's Cheeger and Goresky-MacPherson, independently and
by entirely different methods, discovered that
there is a much larger class of spaces than manifolds that can be
assigned a sequence of ``homology groups'' that satisfy Poincar\'e
duality.  Given  the central role that Poincar\'e duality
plays in surgery theory, it was inevitable that this would lead to
a new environment for the applications of surgery.

\subhead 2.1\endsubhead
Let $X$ be a stratified space where $X^i\setminus X^{i-1}$ is an $i$-dimensional 
$F$-homology manifold, for a field $F$.  We shall assume that the
codimension one stratum is of codimension at least two and that 
$X\setminus X^{n-1}$ (the nonsingular part) is given 
an $F$-orientation; for simplicity
we will also mainly be concerned with the case of $F=\bq$.  It pays to
think PL, as we shall, but see \cite{Q3} for an extension to homotopically
stratified sets.  A
{\it perversity\/} $p$ is a nondecreasing function from the natural numbers to
the nonnegative integers, with $p(1) = p(2) = 0$, and for each $i$, 
$p(i+1)\leq p(i)+1$.  The {\it zero perversity\/} is the identically 
$0$ function and the
{\it total perversity\/} $t$ has $t(i) = i-2$ for $i\geq 2$.
Two perversities, $p$ and $q$ are {\it dual}
if $p+q=t$.

We say that $X$ is {\it normal\/} if the link of any simplex of
codimension larger than $1$ is connected.  This terminology is
borrowed from algebraic geometry.  It is not hard to ``normalize''
``abnormal'' spaces by an analogue of the construction of the
orientation two-fold cover of a manifold.

A {\it chain\/} is just what it always was in singular homology: we
say it is {\it $p$-transverse\/}, or {\it $p$-allowable\/}, 
if for every simplex in the
chain $\Delta\cap X^{n-i}$ has dimension at most $i$ larger than what would be
predicted by transversality and the same is true for the simplices in
its boundary that have nonzero coefficient.

Note:  It is not the case that the chain complex of $p$-transverse
chains with coefficients in $R$ is just that for $\bz$ and then 
$\otimes R$, as it
would be in ordinary homology, because a non-allowable chain can
become allowable after tensoring when some simplex in the
boundary gets a $0$-coefficient.

\subhead 2.2\endsubhead  $IH_p(X)$ is the homology obtained by considering 
$p$-allowable
chains.  It is classical for normal spaces that $IH^t$ is just ordinary
homology; a bit more amusing is the theorem of McCrory that $IH^0$ is
cohomology in the dual dimension.  The forgetful map is capping with
the fundamental class.

Note that $IH$ is {\bf not} set up to be a functor.  It turns out to be
functorial with respect to {\it normally nonsingular\/} or ({\it homotopy\/})
{\it transverse\/} maps.  (We'll discuss these in a great deal more details in
\S\S 4,5.)  Thus, it {\bf is} functorial with respect to (PL)
homeomorphisms and inclusions of open subsets and collared
boundaries.

        Note also, that one can give ``cellular'' versions of $IH$, which
means that one can define perverse finiteness obstructions and
Reidemeister and Whitehead torsions in suitable circumstances.
(See \cite{Cu, Dr}.)  Here one would usually want to build in
refinements to integer coefficients that we will not discuss till
2.10 below.

\subhead 2.3 \endsubhead The main theorems of \cite{GM1} are that 
(1) $IH$ is stratification
independent (indeed it is a topological invariant, even a stratified
homotopy invariant) and (2)  for dual perversities the groups in dual
dimensions are dual.  The latter boils down to Poincar\'e duality in
case $X$ is a manifold, however it is much more general.

\subhead 2.4\endsubhead  
What is important in many applications is that one can often get
a self duality.  Unfortunately, there is no self dual perversity
function (what should $p(3)$ equal?).  However, we have two {\it middle
perversities\/} $0,0,1,1,2,2,\dots$ and $0,0,0,1,1,2\dots$; note that these differ
only on the condition of intersections with odd codimensional strata.
Consequently, for spaces with only even codimensional strata, the
middle intersection homology groups are self dual.

\subhead 2.5\endsubhead  
It turns out that the middle perversity groups have many other
amazing properties.  Cheeger independently discovered the ``De Rham''
version of these.  He gave a polyhedral $X$ as above a piecewise flat
metric (i.e. flat on the simplices, and conelike) and observed that the
$L^2$ cohomology of the incomplete manifold obtained by removing the
singular set was very nice.  Under a condition that easily holds when
one has even codimensional strata, the $\ast$ operator takes $L^2$ forms to
themselves, and one formally obtains Poincar\'e duality.  A
consequence of this is that the  {\bf Kunneth formula holds}.

        In addition, Goresky and MacPherson \cite{GM3} proved that Morse
theory takes a very nice form for stratified spaces when you use
intersection homology.  This leads to a proof of the Lefshetz
hyperplane section theorem.  (A sheaf theoretic proof appears in 
\cite{GM2}.)  \cite{BBD} proved hard Lefshetz 
for the middle perversity intersection
homology of a singular variety using the methods of characteristic $p$
algebraic geometry.  This requires the sheaf theoretic reformulation
to be discussed below.  Finally Saito \cite{Sa} gave an analytic proof of
this and a Hodge decomposition for these groups.

\subhead 2.6\endsubhead  
Let us return to pure topology by way of example.  Consider a
manifold with boundary $W,\partial W$, and the singular space obtained by
attaching a cone to $\partial W$.  Normality would correspond to the
assumption that $\partial W$ is connected.

        What are the intersection homology groups in this case?  Fix $p$.
We would ordinarily not expect any chain of dimension less than $n$ to
go through the cone point.  Once $i+p(i)$ is at least $n$, we begin
allowing all chains to now go through the cone point, so one gets
above that dimension all of the reduced homology.  Below that
dimension, we are insisting that our chains miss the cone point, so
one gets $H_\ast(W)$.  There is just one critical dimension where the chain
can go through and the boundary cannot: here one gets the image of
the ordinary homology in the reduced homology.

        Using these calculations, one can reduce the 
Goresky-MacPherson duality theorem to Poincar\'e-Lefshetz duality for the
manifold with boundary.

        If the dimension of $W$ is even, one gets in the middle dimension
(for the middle perversity) the usual intersection pairing on $(W,\partial W)$
modulo its torsion elements.

        Note though that if $W$ is odd dimensional the failure of self
duality is caused by the middle dimensional homology of $\partial W$.  If its
homology vanished, we'd still get Poincar\'e duality.

\subhead 2.7\endsubhead  
Of course, one immediately realizes that one can now define
signatures for spaces with even codimensional singularities (that
lie in the Witt group $W(F)$ of the ground field.)  We'll, for now, only
pay attention to $F=\br$ and ordinary signature.

        Thom and Milnor's work on PL $L$-classes and Sullivan's work on
$KO[1/2]$ orientations for PL manifolds all just depend on a cobordism
invariant notion of signature that is multiplicative with respect to
products with closed smooth manifolds.  Thus, as observed in \cite{GM I}
it is possible to define such invariants lying in ordinary homology
and $KO[1/2]$ of any space with even codimensional strata.

\subhead 2.8\endsubhead  
It is very natural to sheafify.  Nothing prevents us from
considering the intersection homology of open subsets and one sees
that for each open set one has duality between locally finite
homology and cohomology.  It turns out that the usual algebraic
apparatus of surgery theory mainly requires self dual sheaves rather
than manifolds.  So we can define symmetric signatures that take
the fundamental group into account, which are just the assemblies
(in the sense of assembly maps) of the classes in 2.7.

\subhead 2.9\endsubhead  
The original motivations to sheafify were rather different.
Firstly, using sheaf theory there are simple Eilenberg-Steenrod type
axioms that can be used to characterize $IH$; these are useful for
calculational purposes and for things like identifying the Cheeger
description with the geometric definition of Goresky and
MacPherson.

Secondly, using various constructions in the derived category
of sheaves, push forwards and proper push forwards and truncations
of various sorts, it is possible to give a direct abstract definition of
IH without using chains.  This definition is appropriate to
characteristic $p$ algebraic varieties.

Finally, there is a very simple sheaf theoretic statement,
Verdier duality, that can be used to express locally the self duality
of the intersection homology of all open subsets of a given $X$.  It says
that $IC^m$ is a self-dual sheaf for spaces with even codimensional
singularities.  We will see below that this is quite a powerful
statement.

\subhead 2.10\endsubhead  
We can ask for which spaces is $IC$ self dual?  We know that all
spaces with even codimensional strata have this property, but they
are not all of them, for we saw in 2.6 that if we have an isolated
point of odd codimension one still gets Poincar\'e duality in middle
perversity $IH$ if (and only if) the middle dimensional homology of the
link -- which is a manifold -- vanishes.  One can generalize this
observation to see that if the link of each simplex of odd
codimension in $X$ has vanishing middle $IH^m$, then $IC$ is self dual on $X$.
(Indeed this condition is necessary and sufficient.) Such $X$'s were
christened by Siegel \cite{Si}, {\it Witt spaces\/}.  Actually they were
introduced somewhat earlier by Cheeger as being the set of spaces
for which the $\ast$ operator on $L^2$ forms on the nonsingular part behaves
properly.

The main point of Siegel's thesis, though, was to compute the
bordism of Witt spaces.  Obviously the odd dimensional bordism
groups vanish, because the cone on an odd dimensional Witt space is
a Witt nullcobordism .  For even dimensional Witt spaces this only
works if there is no middle dimensional $IH^m$.  By a surgery process
on middle dimensional cycles, he shows that you can cobord a Witt
space to one of those if 
and only if the quadratic form in middle $IH^m(~ ;\bq)$ is
hyperbolic -- so there is no obstruction in $2$ mod $4$, but there's an
obstruction in $W(\bq)$ in $0$ mod $4$.  Moreover, aside from dimension $0$,
where all that can arise is $\bz\subseteq W(\bq)$ 
given by signature, in all other
multiples of $4$ all the other (exponent $4$ torsion, computed in \cite{MH})
elements can be explicitly constructed, essentially by plumbing.  The
isomorphism of the bordism with $W(\bq)$ is what gives these spaces
their name.

However, making use of the natural transformations discussed
above,  we actually see that Witt spaces form a nice cycle theory for
the (connective) spectrum $L(\bq)$ if we ignore dimension $0$.  Siegel
phrases it by inverting $2$:

\proclaim{Theorem} 
Witt spaces form a cycle theory for connective $KO\otimes\bz[1/2]$,
i.e. 
$$\Omega^{Witt}(X)\otimes\bz[1/2]\to  KO(X)\otimes\bz[1/2]$$
is an isomorphism.
\endproclaim

Pardon, \cite{P} building on earlier work of Goresky and Siegel,
\cite{GS}, showed that the spaces with {\bf integrally self dual $IC$} form a
class of spaces (which does not include all spaces with even
codimensional strata: one needs an extra condition on the torsion of
the one off the middle dimension $IH$ group)  whose cobordism groups
agree with $L^\ast(\bz)$ and then give a cycle theoretic description for the
connective version of this spectrum.

Other interesting bordism calculations for classes of singular
spaces can be found in \cite{GP}.

\subhead 2.11  (Some remarks of Siegel)\endsubhead  
The fact that one has a bordism
invariant signature for Witt spaces contains useful facts about
signature for manifolds.  For instance, using the identification of
signature for manifolds with boundary with the intersection
signature of the associated singular space with an isolated singular
point, it is easy to write down a Witt cobordism (the pinch
cobordism) which proves Novikov's additivity theorem \cite{AS}.

Also, the mapping cylinder of a fiber bundle is not always a
Witt cobordism: there is a condition on the middle homology of the
fiber.  One could have thought that one can still define signature for
singular spaces where the links have signature zero (obviously one
can't introduce a link type with nonzero signature and have a
cobordism invariant signature).  However, Atiyah's example of
nonmultiplicativity of signature gives a counterexample to this \cite{A}.
It is thus quite interesting that having no middle homology {\bf is} enough
for doing this.

\subhead 2.12\endsubhead  
Siegel's theorem has had several interesting applications.  The
first is a purely topological disproof of the integral version of the
Hodge conjecture (already disproven by analytic methods in \cite{AH}) on
the realization of all $(p,p)$ homology classes of a nonsingular variety
by algebraic cycles.  If one were looking for nonsingular cycles, then
one can use failure of Steenrod representability, or better, Steenrod
representability by unitary bordism!, but here we allow singular
cycles.  Thanks to Hironaka, we could apply resolution of
singularities to make the argument work anyway.  However, even
without resolution one sees that these homology classes have a
refinement to $K$-homology, which is a nontrivial homotopical
condition (as in \cite{AH} which develops explicit counterexamples along
these lines).

Another application stems from the fact that the bordism
theory is homology at the prime $2$.  Since one can define a signature
operator for Witt spaces which is bordism invariant \cite{PRW}, one can view
the signature operator from the point of view of Witt bordism and
thus obtain a refinement at the prime $2$ of the $K$-homology class of
the signature operator to ordinary homology \cite{RW}.  This, then
implies that the $K$-homology class of the signature operator is a
homotopy invariant for, say, $RP^n$.

Yet other applications concern ``eta type invariants''.  The basic
idea for these applications is that if one knows the Novikov
conjecture for a group $\pi$, then by Siegel's theorem 
$\Omega^{Witt}(B\pi)\to L(\bq\pi)$ rationally injects.  
This means that one can get geometric
coboundaries from homotopical hypotheses.  Thus, for instance,
homotopy equivalent manifolds should be rationally Witt cobordant
(preserving their fundamental group).

In particular, then, the invariant of Atiyah-Patodi-Singer \cite{APS}
associated to an an odd dimensional manifold with a unitary
representation of its fundamental group can only differ, for
homotopy equivalent manifolds, that a twisted signature of the
cobounding Witt space, e.g. a rational number.  In \cite{W3}, known
results regarding the Novikov conjecture and the deformation
results of \cite{FL} are used to prove this unconditionally.

A similar application is made in \cite{W6} to define 
``higher rho invariants'' for various classes of manifolds.  
For instance, say that
a manifold is antisimple if its chain complex is chain equivalent to
one with $0$ in its middle dimension (this can be detected
homologically).  Then its symmetric signature vanishes and,
therefore, assuming the Novikov conjecture, it is Witt nullcobordant.
By gluing together the Witt nullcobordism and the algebraic
nullcobordism one obtains a closed object one dimension higher,
whose symmetric signature (modulo suitable indeterminacies) is an
interesting invariant of such manifolds.  It can be used to show that
the homeomorphism problem is undecidable even for manifolds
which are given with homotopy equivalences to each other \cite{NW}.

\subhead 2.13  (Dedicated to the Cheshire cat)\endsubhead  
Associated to any Witt space
one has a self dual sheaf, namely $IC^m$.   Actually, the cobordism
group of self dual sheaves over a space $X$ (assuming the self duality
is symmetric) can be identified with $H_\ast(X;L^\ast(\bq))$, 
(see \cite{CSW} for a
sketch, and \cite{Ht} for a more general statement including some more
general rings\footnote{In general there are algebraic 
$K$-theoretic difficulties with identifying the
Witt group self dual sheaves, at $2$, with a homology theory.  However, as Hutt
himself was aware, one can certainly include many more rings than included
in that paper.}.

This statement has some immediate implications:  Since $IC^m$
is topologically invariant, all of the characteristic classes
introduced for singular spaces in this way are topologically
invariant.  (This is basically the topological invariance of rational
Pontrjagin classes extended to Witt spaces.)

We thus have, away from $2$, three rather different descriptions
of $K$-homology:  Witt space bordism, homotopy classes of abstract
elliptic operators (see \cite{BDF, K}), and bordism of self dual sheaves
(and, not so different from this one: controlled surgery obstruction
groups).

A number of applications to equivariant and stratified surgery
come from these identifications (and generalizations of them).  We
will return to some of these in \S6.

\subhead 2.14\endsubhead  
A very nice application of cobordism of the self dual sheaves
associated to IH and its various pushforwards is given in \cite{CS2}.  The
goal is to extend the usual multiplicativity of signature in fiber
bundles (with no monodromy) to stratified maps.  We will not give a
precise definition of a stratified map, but it is the intuitive notion,
e.g. a fiber bundle has just one stratum.

\proclaim{Theorem}  Let $f:X\to Y$ be a stratified map betweens spaces with
even codimensional strata, and suppose  that all the strata of $f$ 
are
of even codimension and the pure strata are simply connected.  We
then have
$$f_\ast(L(X)) = \sum sign(c(star_f(V)))L(V) $$
where $V$ runs over the strata of $f$ (which is a substratification of $Y$).
Here $c(\cdot)$ stands for a compactification -- in this case it means the
following.  If $V = Y$ it is just the generic fiber.  If $V$ is a proper
stratum, then one can consider $f^{-1}(cone(L))$, where $L$ is the link of a
generic top simplex of $V$, and then one point compactify it (=cone off
its boundary).
\endproclaim

One can deal with nonsimply connected open strata by putting
in a correction term for the monodromy action of 
$\pi_1(\inr(V))$ on $IH(c(star_f(V)))$.

The proof of the theorem in \cite{CS2} is very pretty; it makes use
of the machinery on perverse sheaves found in \cite{BBD} but in explicit
cases essentially produces an explicit cobordism between $f_\ast IC(X)$
and an explicit sum of other intersection sheaves: one for each
stratum of the map.

\remark{Remarks}
\roster  
\item  To get a feeling for the theorem it is worth considering
a few special cases.  Firstly, the case of a fiber bundle just reduces
to \cite{CHS}.  As a second special case, if one considers a pinch map
from a union of two manifolds along their common boundary, the
formula boils down to Novikov additivity, and the cobordism implicit
in the proof is the pinch cobordism of 2.11.  As a final example, one
can consider the case of a circle action on a manifold.  Aside from
some fundamental group points, there is a similar cobordism
between M and some projective space bundles over the fixed set
components of the circle action, and the formula of the theorem
generalizes by considering projection to the quotient -- with some
slight modifications for $0$ mod $4$ components of fixed set, which
lead to non-Witt singularities -- (or specializes to) the formula in
\cite{W2} that identifies the higher signature of a manifold and that of
its fixed point set of any circle action with nonempty fixed point
set.  The cobordism (discussed in both \cite{W2} and \cite{CS2}) is then the
bubble quotient cobordism of \cite{CW3}.
\item  In the case of an {\bf algebraic map}, one could directly apply 
\cite{BBD}
which gives a beautiful and deep decomposition theorem for $f_\ast IC(X)$
and the general machinery on self dual sheaves to prove the
existence of a formula like the one in the theorem.  However, it is
not so clear what the coefficients are.
\item  Still in the algebraic case, it is important to realize that there
are many additional characteristic classes that can be defined for
singular varieties beyond just the $L$-classes, for instance,
MacPherson Chern classes and Todd classes.  In \cite{CS3}, there are
announced generalizations of the basic formula where the meaning
of c is different:  one must use projective completion to get a
variety, and then the formula must be rewritten to account for the
extra topology (think about the case of intersection Euler
characteristic classes, which can be dealt with by the proof of the
usual Hurwitz formula for Euler characteristic of branched cover,
together with the sheaf version of intersection Kunneth).  To prove
such formulae one uses deformation to the normal cone (see \cite{Fu2}) to
replace the cobordism theory.
\endroster
\endremark

\subhead 2.15\endsubhead  
It is worth mentioning but beyond the purview of this survey to
describe in any detail the applications of 2.14 given in \cite{CS3, CS4, Sh2} to
lattice point problems.  The connection goes via the theory of toric
varieties for which there are several excellent surveys \cite{Od, 
Da, Fu1}, which gives an assignment of a (perhaps singular)
variety to every convex polygon with lattice point vertices on
which a complex torus acts.  (See also \cite{Gu} for a discussion of the
purely symplectic aspects of this situation.)  Problems like counting
numbers of lattice points inside such a polytope (= computation of
the Erhart polynomial) and Euler MacLaurin summation formulae can
be reduced to calculations of the Todd class, which are studied in
tandem with $L$-classes using the projection formulae.  These, in
turn, have substantial number theoretic implications.

\head 3. The geometry of homotopically stratified spaces
\endhead
One of the strengths of Quinn's formulation of manifold homotopically
stratified spaces is that the defining conditions are homotopic theoretical
(except, obviously, the geometric condition that the strata be manifolds).
This, of course, makes it easier to verify the conditions, but harder
to establish geometric facts about manifold homotopically stratified spaces.
Nevertheless, Quinn was able to prove two
important geometric results:  homogeneity and stratified
$h$-cobordism theorems. 

Quinn's homogeneity result says that if $x,y$ are two points in the same
component of a stratum (with adjacent strata of dimension at least $5$)
of a manifold homotopically stratified space $X$, then there is a 
self-homeomorphism (in fact, isotopy)
of $X$ carrying $x$ to $y$. Quinn obtains this as a consequence of an
stratified isotopy extension theorem (an isotopy on a closed union of
strata can be extended to a stratum preserving isotopy on the whole space).
In turn, Quinn proves the isotopy extension theorem by using the
full force of his work on ``Ends of maps'' (see \cite{Q2,IV}).

As an example of the usefulness of the homogeneity result,
consider a finite group acting on a manifold $M$. Even though the action
need not be locally linear, the quotient $M/G$ is often a manifold
homotopically stratified space. 
Thus, the homogeneity result can be used to verify local linearity
by establishing local linearity at a single point of each stratum
component.
Quinn first used this technique to construct locally linear actions
whose fixed point set does not have an equivariant mapping cylinder
neighborhood \cite{Q2,II}. Weinberger \cite{W1}
and Buchdahl, Kwasik and Schultz \cite{BKS} have also used this result
to verify that certain actions were locally linear.

It turns out that there is an alternative way to prove Quinn's homogeneity
theorem which is based on engulfing (the classical
way that homotopy information is converted into homeomorphism 
information in manifolds).
In fact, this alternative method uses a description of neighborhoods of
strata together with  MAF (manifold approximate fibration) technology, and is
useful for many other geometric results.

We have seen in \S1 that in certain formulations of conditions on a
stratification $\Sigma=\{ X_i\}$ of a space $X$ one considers
tubular maps
$$\tau_i:U_i\to X_i\times [0,+\infty)$$
where $U_i$  is a neighborhood of $X_i$ and
$\tau_i$ restricts to the identity $U_i\setminus X_i\to X_i\times (0,+\infty)$.
For Whitney stratifications, the tubular maps are submersions on each 
stratum and fibre bundles over $X\times (0,+\infty)$. Since strata of
manifold homotopically stratified spaces need not have mapping cylinder
neighborhoods, such a result is too much to hope for in general.
However, there is the 
following result which was proved by Hughes, Taylor, Weinberger and
Williams \cite{HTWW} in the case of two strata and by Hughes \cite{H3}  
in general.

\proclaim{Theorem} For manifold homotopically stratified spaces in which all
strata have dimension greater than or equal to $4$, tubular maps
exist which are manifold stratified approximate fibrations.
\endproclaim

The neighborhoods of the strata which are the domains of these MSAF (manifold
stratified approximate fibration) tubular maps are called
{\it teardrop neighborhoods\/}. They are effective substitutes for
mapping cylinder neighborhoods, and the result should be thought of as
a {\it tubular neighborhood theorem\/} for stratified spaces. 
The point is  that even though Quinn's definition does not postulate 
neighborhoods with any kind of reasonable tubular maps, one is able to
derive their existence. The situation is optimal: minimal conditions in
the definition with much stronger conditions as a consequence.
This makes the surgery theory conceptually easier than for
geometrically stratified sets for which the geometric neighborhood
structure must be part of the data.

As mentioned above, these teardrop neighborhoods can be used to 
give a different proof of Quinn's isotopy extension theorem.
Manifold approximate fibrations have the approximate isotopy covering
property \cite{H1}. This property holds in the stratified setting
and is used inductively to extend isotopies from strata to their
teardrop neighborhoods. In fact, parametric isotopy extension is now
possible whereas Quinn's methods only work for a single isotopy.

Similarly, other results in geometric topology can be extended to
manifold homotopically stratified sets by using MAF technology.
For example, the homeomorphism
group of a manifold homotopically stratified set is locally contractible,
and a stratified version of the Chapman and Ferry \cite{ChF} 
$\alpha$-approximation theorem holds.
In short, the case can be made that manifold homotopically stratified
sets are the correct
topological version of polyhedra and Thom's stratified sets.

\head 4. Browder-Quinn theory
\endhead
In \cite{BQ}, Browder and Quinn introduced an interesting, elegant, and 
useful general classification theory for strongly stratified spaces.
The setting is a category where one has a fixed choice of strong 
stratification as part of the data one is interested in.

\subhead 4.1
\endsubhead
Let $X$ be a strongly stratified space (e.g. a geometrically stratified
space as in \S1.3) with {\it closed pure
strata} $\ov  X^i$ (see \S1.3). 
An $h$-cobordism with boundary $X$ is a stratified
space $Z$ with boundary $X\cup X'$ where the inclusions of $X$ and
$X'$ are stratified homotopy equivalences, and the neighborhood data
for the strata of $Z$ are the pullbacks with respect to the retractions of
the data for $X$ (and of $X'$). This condition is automatic in the
PL and Diff categories when one is dealing with something like the
quotient of a group action stratified by orbit types.

\proclaim{Theorem} The $h$-cobordisms with boundary $X$ \rom(ignoring
low dimensional strata\rom) are in a 1--1  correspondence with a group 
$\Wh^{BQ}(X)$. There is an isomorphism 
$\Wh^{BQ}(X) \cong \oplus \Wh(\ov X^i)$.
\endproclaim

The map $\Wh^{BQ}(X) \to \oplus \Wh(\ov X^i)$ is given by sending
$(Z,X)\to (\tau (\ov Z^i,\ov X^i))$.
One proves the isomorphism (and theorem) inductively, using the
classical $h$-co\-bord\-ism theorem to begin the induction, and using
the strong stratifications to pull up product structures to
deal with one more stratum.

\subhead 4.2
\endsubhead
The surgery theory of Browder and Quinn deals with the problem of turning
a degree one normal map into a stratified homotopy equivalence which is
{\bf transverse}, i.e. one for which the strong
stratification data in the domain is the pullback of the data from the
range.

This transversality is, in practice, the fly in the ointment. When one is
interested in classifying embeddings or group actions usually the bundle data
is something one is interested in understanding rather than {\it a priori\/}
assuming. Still, in some problems (e.g. those mentioned in 6.3) one can 
sometimes prove that transversality is automatic. Also, of course, if one uses 
the machinery to construct examples, it is certainly fine if one produces
examples that have extra restrictions on the bundle data.

\subhead 4.3
\endsubhead
Either by induction or by mimicking the usual identification of normal
invariants, one can prove that 
$NI^{BQ}(X)\cong [X;F/Cat]$.

\subhead 4.4 
\endsubhead
Define $S^{BQ}(X)$ to be the strongly stratified spaces with a
transverse stratified simple homotopy equivalence to $X$ up to 
$Cat$-strongly stratified simple isomorphism (note this implicitly is keeping
track of ``framings''). Then, one has groups $L^{BQ}(X)$ and a long exact
surgery sequence:
$$\cdots \to L^{BQ}(X\times I\/\rel\partial )\to S^{BQ}(X)\to [X;F/Cat]\to
L^{BQ}(X).$$

\subhead 4.5
\endsubhead
Note that unlike the Whitehead theory $L^{BQ}(X)$ does not decompose into a
sum of $L$-groups of the closed strata. Indeed, for a manifold with
boundary $S^{BQ}$ is just the usual structure set
(existence and uniqueness of collars gives the strong stratification
structures) and the $L$-group is the usual $L$-group of a manifold
with boundary, i.e. is a relative $L$-group, not a sum of absolute groups.

However, there is a connection between the $L$-groups of the pure strata
and $L^{BQ}(X)$. This exact sequence generalizes the exact sequence of a 
pair in usual $L$-theory, and expresses the fact that as a space
$L^{BQ}(X)$ is the fiber of the composition
$$L(X_0)\to L(\partial\text{Neighborhood}(X_0))\to L(\cl(X\setminus X_0))$$
where the first map is a transfer and the second is an inclusion.

\subhead 4.6 
\endsubhead
The proof of this theorem is by the method of chapter 9 of 
\cite{Wa}: one need only verify the $\pi-\pi$ theorem. This is done by
induction.

\head 5. Homotopically stratified theory
\endhead
If one does not want to insist on the transversality condition
required in the Browder-Quinn theory,
or if one is only dealing with homotopically stratified spaces, it is
necessary to proceed somewhat differently.  For more complete
explanations, see \cite{W4}, \cite{W5}.  We will only discuss the topological
version.  The PL version is simpler but slightly more complicated.

\subhead 5.1  The $h$-cobordism theorem\endsubhead
That new phenomena would arise in any systematic study of
Whitehead torsion on nonmanifolds was clear from the start.
Milnor's original counterexamples to the Hauptvermutung for
polyhedra were based on torsion considerations \cite{M1}.  Siebenmann
gave examples of locally triangulable nontriangulable spaces -- not
at all due to Kirby-Siebenmann considerations, but rather $K_0$.  More
pieces came forward in the work on Anderson-Hsiang \cite{AnH1, AnH2} and
then in \cite{Q2}, which showed that under appropriate K-theoretic
hypotheses one can triangulate, and therefore apply the
straightforward PL torsion theory.

        Real impetus came from the theory of group actions.  Cappell
and Shaneson \cite{CS1} gave the striking examples of equivariantly
homeomorphic representation spaces, which laid down the gauntlet
to the topological community at large to deal with the issue of
equivariant classification.  $h$-cobordism theorems suitable for the
equivariant category were produced by Steinberger (building on joint
work with West) \cite{St} and by Quinn \cite{Q4} in the generality of
homotopically stratified spaces (although the theorem in that paper
does not include realization of all torsions in an $h$-cobordism
\footnote{As pointed out in \cite{HTWW}, 
the teardrop neighborhood theorem can be used to
complete the proof of realization.}).

        The ultimate theorem asserts, as usual, that (ignoring low
dimensional issues)  $h$-cobordisms on a stratified space X are
classified by an abelian group $\Wh^{top}(X)$.

\proclaim{Theorem}  $\Wh^{top}(X)\cong\oplus \Wh^{top}(X^i,X^{i-1})$, 
and we have an exact
sequence
$$\multline\cdots\to H_0(X^{i-1}; \text{\bf Wh}(\pi_1(\ho))
\to \Wh(\pi_1(X^i\setminus X^{i-1}))
\to\\
\Wh^{top}(X^i,X^{i-1})\to H_0(X^{i-1}; {\text{\bf K}}_0(\pi_1(\ho))\to 
K_0(\pi_1(X^i\setminus X^{i-1}))
\endmultline$$
\endproclaim

Boldface terms are spectra.  This decomposition of $\Wh^{top}$ into a
direct sum does not respect the involution obtained by turning 
$h$-cobordisms upside down, is a pleasant descendent of the analogous
fact in the Browder-Quinn theory.  It does not have an analogue in 
L-theory.

\subhead 5.2  Stable classification\endsubhead
Ranicki (following a sketch using geometric Poincar\'e
complexes in place of algebraic ones, by Quinn) has elegantly
reformulated the usual Browder-Novikov-Sullivan-Wall surgery
exact sequence in the topological manifold setting as the assertion
that there is a fibration:
$$\text{\bf S}(M)\to \text{\bf H}_n(M;\text{\bf L}(e))\to
\text{\bf L}_n(\pi_1(M))$$
where {\bf X} means a space (or better a spectrum) whose homotopy
groups are those of the group valued functor ordinarily denoted by $X$.
$S(M)$ is the {\it structure set\/} of $M$, 
which essentially\footnote{In actuality, for our purposes it is best 
to use the finite dimensional ANR
homology manifolds, and the equivalence relation is $s$-cobordism.  See
Mio's paper \cite{Mi}
in this volume for a discussion of the difference this makes.  (It is at most a
single $\bz$ if $M$ is connected.)} 
consists of the
manifolds simple homotopy equivalent to $M$ up to homeomorphism.
The map  $\text{\bf H}_n(M;\text{\bf L}(e))\to\text{\bf L}_n(\pi_1(M))$ 
is called the {\it assembly map\/} and
can be defined purely algebraically.  Geometrically it has several
interpretations: most notably, as the map from normal invariants to
surgery obstructions in the topological category, or as a forgetful
map from some variant of controlled surgery to uncontrolled
surgery.

Since the assembly map has a purely algebraic definition, one
can ask whether it computes anything interesting if $X$ is not a
manifold? and alternatively, if $X$ is just a stratified space, what is
the analogue of this sequence?

Cappell and the second author gave an answer to the first
question in \cite{CW2} where it is shown (under some what more
restrictive hypotheses) if $X$ is a manifold with singularities, i.e. $X$
contains a subset $\Sigma$ whose complement is a manifold, and suppose
further that $\Sigma$ is $1$-LCC embedded\footnote{This 
means that maps of $2$-complexes into $X$ can be 
deformed slightly to miss $\Sigma$.}  in $X$, 
then $S^{alg}(X)\cong S(X \rel\Sigma)$ where
$S^{alg}(X)$ denotes the fiber of the algebraically defined assembly map
$H_\ast(X;L(e))\to L(X)$ and
$S(X \rel\Sigma)$ means 
$$\split\{\varphi:X'\to X ~|&~ \varphi \text{ is 
a stratified 
simple\footnotemark [5]\ homotopy
equivalence}\\ 
&\text{ with $\varphi|\Sigma(X')$ already a homeomorphism}\}.
\endsplit$$
\footnotetext [5]{The materiel of 5.1 can be used to make sense of this.}%
The answer to the second question is a bit more complicated,
and actually requires two fibrations in general.  The first  is a
stable generalization of the surgery exact sequence:
$$S^{-\infty}(X \rel Y)\to\text{\bf H}_0(X;\text{\bf L}^{BQ-\infty}
(--,\rel Y))\to\text{\bf L}^{BQ-\infty}(X \rel Y).$$
Here the superscripts $-\infty$ denote that we are using a stable version
of structure theory: we will soon explain that it only differs from
the usual thing at the prime $2$, and the phenomenon is governed by
algebraic K-theory.  The coefficients of the homology is with
respect to a cosheaf of L-spectra: to each open set $U$ of $X$ one
assigns the spectrum $\text{\bf L}(U \rel U\cap Y \text{with compact support})$.
The $BQ$
superscripts are a slight generalization of the theory discussed in 
\S4.

To complete the theory one needs a destabilization sequence.
This is given by the following:
$$S(X) \to S^{-\infty}(X)\to \hat{H}(\bz_2;\Wh^{top}(X)^{\leq 1})$$
Here $S(X)$ is the geometric structure set, and $S^{-\infty}(X)$ is the
stabilized version, which differ by a Tate cohomology term.  An
analogous sequence developed for a quite similar purpose appears in
\cite{WW}.  Indeed in \cite{HTWW} the theory outlined in this subsection is
deduced from the \cite{WW} results using blocked surgery \cite{Q1, BLR, CW2}
and \cite{HTW1,2} (the classification of manifold approximate fibrations)
and the teardrop neighborhood theorem.  On the other hand, there are
different approaches to all this using controlled end and/or surgery
theorems that are sketched in \cite{W4}, \cite{W5}.

\head 6.  Some applications of the stratified surgery exact sequence
\endhead
In practice the application of the theory of the previous 
section requires additional input for the calculations to be either 
possible or comprehensible.  See \cite{CW4} for the application to topological 
group actions.  The last 100 pages of \cite{W4} also gives more 
applications than we can hope to discuss here.

\subhead 6.1 
\endsubhead
Probably the simplest interesting and illustrative example of the
classification theorem is to locally flat topological embeddings.  The
first point is that this problem is susceptible to study by these methods:
Every topological locally flat embedding gives a two stratum homotopically
stratified space where the holink is a homotopy sphere {\it and conversely\/}.
This last is essentially Quinn's characterization of local flatness in 
\cite{Q2,I}.

Things are very different in codimensions one and two from
codimension three and higher.  We will defer to 6.3 the low codimensional
discussion and restrict our attention here to the last of these cases.

\proclaim{Lemma}  If $(W,M)$ is a manifold pair with $cod(M)>2$, then 
$L^{BQ}(W,M)\cong
L(W)\times L(M)$ where the map is the forgetful map.
\endproclaim

The proof is straightforward.  Note that the lemma implies the
analogous statement holds at the level of cosheaves of spectra ($\cong$ being
quasi-isomorphism).  It is quite straightforward in this case to work out the
Whitehead theory:  there are no surprises.  Thus, we obtain:

\proclaim{Corollary}  $S(W,M)\cong S(W)\times S(M)$.
\endproclaim

Note that the discussion makes perfect sense even if $(W,M)$ is just
a Poincar\'e pair (see \cite{Wa}), and then the corollary boils down to the
statement that isotopy classes of embeddings of one topological manifold in
another (in codimension at least $3$) are in a $1-1$ correspondence with the
Poincar\'e embeddings (see \cite{Wa}).

(Actually, a bit more work enables one to prove the same thing for
$M$ an ANR homology manifold.)

\subhead 6.2 \endsubhead 
Using the material from \S2 we can also analyze, away from $2$, $S(X)$ for
a very interesting class of spaces that have even codimensional strata.  We
continue to let $S^{alg}(X)$ denote 
the fiber of the classical assembly map 
$H_\ast(X;L(e))\to L(X)$.
It is what the structure set of $X$ would be if $X$ were a manifold.

\proclaim{Theorem(\cite{CW2})}  
If $X$ is a space with even codimensional strata and all
holinks of all strata in one another simply connected, then there is an
isomorphism $[1/2]$
$$S(X)\cong\oplus S^{alg}(V)$$
where the sum is over closed strata.
\endproclaim

The proof consists of building an isomorphism $L^{BQ}(X)\cong\oplus L(V)
[1/2]$ for arbitrary $X$ satisfying the hypotheses.  It is obvious enough how
to introduce $\bq$ coefficients into $L^{BQ}$.  
Ranicki \cite{R} has shown that
introducing coefficients does not change $L$ at the odd primes, but with
$\bq$-coefficients one can make forgetful maps to $L(V;\bq)$ 
using the intersection
chain complexes.

\subhead 6.3\endsubhead  
To give an example where things work out differently, we shall assume
that the holinks are aspherical, and that the Borel conjecture holds for
the fundamental groups of these holinks. (This example is a special case of
the theory of ``crigid holinks'' in \cite{W4}.)

In this case there is nothing good that can be said about the
global $L^{BQ}$ term, in general.  However, the assumptions are enough to imply
that $H_\ast(X;L^{BQ})\cong [X; L(e)]$.  (See \cite{W4}, \cite{BL} for a 
discussion.)  In
particular, for locally flat embeddings in codimension $1$ and $2$, one sees
that the fiber of $S(W,M)\to S(W)$ only involves fundamental group data,
not, say, the whole homology of the manifold and submanifolds.  This, too,
reflects phenomena already found in Wall's book \cite{Wa}.

Another amusing example is $X=$ simplicial complex, stratified by
its triangulation.  Then one gets $L^{BQ}(X)\cong [X; L(e)]$.

There are other interesting examples that display similar phenomena
that come up from toric varieties. The theory of multiaxial actions (see \S2 and
\cite{Dv2}) is another situation where the 
local cosheaves tend to decompose into pieces that can be written in 
simple terms involving skyscraper $L(e)$-cosheaves.  Not all holinks 
are crigid and consequently different phenomena appear: indeed 
signatures 0 and 1 alternate in the simply connected holinks, with 
quite interesting implications.  As a simple exercise, one can see 
that while $S^{6n-1}/U(n)$ is contractible, its structure set
\footnote [6]{Actually, the structure set is $\bz\times\bz_2$ 
with the extra $\bz$ corresponding to actions on 
nonresolvable homology manifolds that are homotopy spheres.}
$S(S^{6n-1}/U(n))\cong\bz_2$.  Similarly, $S^{12n-1}/Sp(n)$ is contractible, 
but its structure set
\footnote [7]{Same caveat as above.} 
is $\bz$. 

\remark{Remark}  If all holinks are simply connected (as in the case of 
multiaxial actions of $U(n)$ and $Sp(n)$) one always has a spectral 
sequence computing $S(X)$ in terms of the $S^{alg}(X_i)$.  For instance, if 
there are just two strata $X\supset\Sigma$, there is an exact sequence: 
$$\dots\to S^{alg}(\Sigma\times I)\to S^{alg}(X)\to S(X)\to S^{alg}(\Sigma)
\to\dots$$
The sequence continues to the left in the most obvious way.  On the 
right it continues via deloops of the algebraic structure spaces.  The 
map $S^{alg}(\Sigma\times I)\to S^{alg}(X)$ 
depends on the symmetric signature of the 
holink (and on the monodromy of the holink fibration).  The case 
where the simply connected holink is rigid is essentially that of 
manifolds with boundary.  The normal invariant term here is $[X; L(e)]$, 
but thought of here as $H(X,\partial X; L(e))$.  On the other hand, in 6.2 we 
gave an important case where this spectral sequence degenerates 
(at least away from the prime 2).
\endremark

\subhead 6.4\endsubhead
  As our final example, let us work out in detail a case that is 
somewhat opposite to the one of the previous paragraph:  $X =$ the 
mapping cylinder of even a PL block bundle $V\to N$, with fiber $F$, 
where $N$ is a sufficiently good aspherical manifold.  (Sufficiently 
good is a function of the reader's knowledge.  Even the circle is a 
case not devoid of interest.)
We are interested in understanding what the general theory 
tells us about $S(X \rel V)$.  

Firstly, there is the calculation of the Whitehead group.  (Or 
even pseudoisotopy spaces....).  In this case, the sequence boils down 
to:
$$ H_0(N;\text{\bf Wh}_1(F))\to \Wh_1(V)\to \Wh^{top}(X \rel V)
\to H_0(N;\text{\bf K}_0(F))
\to K_0(V)$$
In a totally ideal world, the assembly maps $H_0(N;\text{\bf Wh}_1(F))
\to \Wh_1(V)$ 
and $H_0(N;\text{\bf K}_0(F))\to K_0(V)$ would be isomorphisms, and 
$\Wh^{top}(X \rel V)$
would vanish.  However, even in the case of $N = S^1$ where the 
assembly map (for the product bundle) was completely analyzed by 
\cite{BHS}, this is not true.  In that case, there is an extra piece called 
$Nil$ that obstructs this; however, $Nil$ is a split summand.  Thus, the 
assembly maps are still injections, and one obtains an isomorphism 
of $\Wh^{top}(X \rel V)$ with a sum of $Nil$s.  In general, the pattern 
discovered by Farrell and Jones \cite{FJ} shows that the cokernel of 
these assembly maps is at least reasonably conjectured to be a 
``sum'' of $Nil$s.  

The splitting of the $K$-theory assembly map essentially boils 
down to the assertion that $\Wh^{BQ}(X \rel V)\to \Wh^{top}(X \rel V)$ has a 
section.  There are fairly direct proofs of this fact when N is a 
nonpositively curved Riemannian manifold in \cite{FW} and in \cite{HTW3}.  
The first approach notes that putting a PL structure on stratified 
spaces can be viewed (essentially following \cite{AnH1, AnH2}) as a problem of 
reducing the the tangent microbundle to the group of block bundle 
maps:  but in the presence of curvature assumptions this can be done 
in the large by the methods of controlled topology.  

The approach in \cite{HTW3} depends on the same controlled 
topology, but its focus is showing that one can associate a MAF 
structure to any map whose homotopy fiber is finitely dominated.  
The teardrop neighborhood theorem of course provides the relation 
between these approaches.

The same analyses can be done for the (stable) structure set 
$S(X \rel V)$.  In this case one does often have the vanishing of the 
analogue of $Nil$ (although if there's orientation reversal or 
complicated monodromy in the bundle, this might not be the case).  
The structure set is here described as the fiber of the assembly 
map, and thus it often vanishes.

This has an interesting interpretation.  Let us suppose that the 
fiber is $K$-flat, i.e. that $\Wh(\pi_1(F)\times\bz^k) = 0$ 
for all $k$ to avoid any 
potential end obstructions.  In this case one also knows that all 
MAF's are equivalent to block bundle projections.

The vanishing of $S(X \rel V)$ means that $S(X)\cong S(V)$ by the 
``obvious'' fibration:  $S(X \rel V)\to S(X)\to S(V)$.  
(We'll discuss the ``$~$''
marks in a moment.)  Now $S(X)$ is basically the same thing as the 
$F$-block bundles on $N$ with fiber a manifold homotopy equivalent
to $F$.  Thus we have a 
generalized fibration theorem for manifolds with maps to $N$.  
(Indeed, 
the Farrell fibration theorem \cite{Fa} is all that is necessary to feed 
into the machinery to get out the calculation of $L$-groups:  that's the 
content of Shaneson's thesis \cite{Sh1}!)

Without the $K$-flatness, we see that there are still only $Nil$ 
obstructions to obtaining MAF structures (but genuine $K$-theory 
obstructions to getting block structures).

To return to the ``obvious'' fibration, a little thought shows that it is not 
at all obvious.  What is obvious is that it is a fibration over the 
components of $S(V)$ in the image of the map $S(X)\to S(V)$.  We are 
asserting, after the arguments given above, that this image is all 
the components, but prima facie, the argument in whole is circular.  

However, that is not the case as a consequence of the complete 
general theory.  The map $S(X)\to S(V)$ is actually an infinite loop 
map, isomorphic to its own $4^{th}$ loop space (see \cite{CW1, We5}).  Thus, 
the fact that we knew exactness at the $\pi_i$ level for $i=3,4$ gives us 
everything we want for any such ad hoc component problem.  (This is 
exactly the same point involved in continuing the exact sequence of 
6.3 further to the right.)

%%%%%%%%%%%%%%%%%%%%%%%%%%%%%%%%%%%%%%%%%%%%%%%%%%%%%%%%%%%%%%%%%%%%%%%%%%%%%%
%% This is the References section for Hughes and Weinberger
%% ``Surgery and stratified spaces'' for the Wall volume.
%%
%%
%%
%AMS-TEX style master bibliography. 

\enddocument